\documentclass[preprint,12pt]{elsarticle}

\usepackage{amssymb}
\usepackage{amsmath}

\usepackage[T1]{fontenc}
\usepackage[utf8]{inputenc}

\usepackage{amsthm}
\usepackage{xcolor}

\usepackage{listings}
\usepackage{algorithm}
\usepackage[noend]{algpseudocode}

\usepackage{siunitx}

\usepackage{tikz}
\usepackage{tcolorbox}

\usepackage{booktabs}

\usepackage{hyperref}

\journal{Computers and Fluids}

\begin{document}

\begin{frontmatter}



\title{Bathymetry reconstruction from experimental data using PDE-constrained optimisation}

\author[label1]{Judith Angel\corref{cor}}
\author[label2]{J\"orn Behrens}
\author[label1]{Sebastian G\"otschel}
\author[label3]{Marten Hollm}
\author[label1]{Daniel Ruprecht}
\author[label3]{Robert Seifried}

\affiliation[label1]{organization={Chair Computational Mathematics, Institute of Mathematics, Hamburg University of Technology},
            addressline={Am Schwarzenberg-Campus 3},
            postcode={Hamburg 21073},
            country={Germany}}

\affiliation[label2]{organization={Department of Mathematics, Universit\"at Hamburg},
            addressline={Bundesstra{\ss}e 55},
            postcode={Hamburg 20146},
            country={Germany}}

\affiliation[label3]{organization={Institute of Mechanics and Ocean Engineering, Hamburg University of Technology},
            addressline={Ei{\ss}endorfer Stra{\ss}e 42},
            postcode={Hamburg 21073},
            country={Germany}}

\cortext[cor]{Corresponding author. \textit{E-mail address: }\url{judith.angel@tuhh.de}}

%

\begin{abstract}
Knowledge of the bottom topography, also called bathymetry, of rivers, seas or the ocean is important for many areas of maritime science and civil engineering.
While direct measurements are possible, they are time consuming and expensive.
Therefore, many approaches have been proposed how to infer the bathymetry from measurements of surface waves.
Mathematically, this is an inverse problem where an unknown system state needs to be reconstructed from observations with a suitable model for the flow as constraint.
In many cases, the shallow water equations can be used to describe the flow.
While theoretical studies of the efficacy of such a PDE-constrained optimisation approach for bathymetry reconstruction exist, there seem to be few publications that study its application to data obtained from real-world measurements.
This paper shows that the approach can, at least qualitatively, reconstruct a Gaussian-shaped bathymetry in a wave flume from measurements of the water height at up to three points.
Achieved normalized root mean square errors (NRMSE) are in line with other approaches.
\end{abstract}

\begin{graphicalabstract}
\end{graphicalabstract}

\begin{highlights}
\item We use PDE-constrained optimisation to reconstruct a Gaussian-shaped bathymetry in a wave flume
\item Point measurements from two sensors are enough to produce a qualitatively correct reconstruction
\end{highlights}


\end{frontmatter}


\section{Introduction}\label{sec:intro}
Knowing topography of river beds or sea or ocean floors, also called bathymetry, is important for many areas of science and engineering.
Two of many examples would be the prediction of meteorological tsunamis~\cite{Rabinovich2020} or the prediction of the circulation of hazardous materials in estuaries~\cite{Sanders2000}. 
Knowledge of shallow water bathymetry is also important for nautical navigation to avoid stranding.

Various methods exist to measure bathymetry. 
Using sonar measurements from research vessels like ships or remotely operated vehicles, the bathymetry can be measured on a width as much as twice the water depth~\cite{smith2004conventional}. 
However, this technique is time-consuming, expensive and there is the risk of stranding~\cite{sagawa2019satellite}. 
In a second technique, the bathymetry is measured by an airborne LIDAR (Light Detection And Ranging) system. 
Airborne LIDAR bathymetry systems can be used to remotely collect accurate, high-density measurements for both bathymetry and topography in coastal regions~\cite{irish1999scanning}. 
However, this technique is still expensive and only works in areas that allow flying~\cite{sagawa2019satellite}. 
As an alternative approach, satellite images can be used to obtain information about the bathymetry of areas that are difficult to reach by airplane or boat. 
Depth information, which are obtained by analyzing satellite images, is called satellite derived bathymetry (SDB)~\cite{ashphaq2021review}. 
Using satellites, bathymetry information of a large area can be tracked much faster than tracked by oceanographic vessels~\cite{smith2004conventional}. 
However, data quality and pre-processing steps comprise a critical part of the SDB. 
Special attention is needed, for example, regarding the atmospheric effects, sun reflection, and the ripple of the sea surface~\cite{evagorou2022evaluation}. 

Bathymetry information can also be obtained from measurements of the height of waves on the water surface.
The advantage of this approach is that it only requires sea surface measurements which are often easier and cheaper to obtain than measurements of the sea floor.
Mathematically speaking, finding the bathymetry from measuring surface waves is an inverse problem where an unknown state of a system is being inferred from observations.
There are a variety of approaches in the literature how to do that.
Many papers use a form of minimisation with forward and adjoint models to compute a descent direction.
Beckers et al.~\cite{BeckersEtAl2019} use dual-weighted residual method with artificial diffusion for advection and advection-diffusion equation discretised with a discontinuous Galerkin (DG) method to smoothen the forward as well as the adjoint solution.
Khan et al.~\cite{Khan2021,Khan2022} reconstruct the bathymetry using the one-dimensional shallow water equations (SWE) and variational assimilation. 
They show that using more observation locations improves the results and that the application of a low-pass filter on the algorithm removes small-scale noise. Furthermore, they analyse the sensitivity of surface waves to errors in the bathymetry reconstruction.
Sanders and Katopodes~\cite{Sanders2000} derive a continuous adjoint problem from the conservation form of the shallow water equations with nonlinear bottom friction term and the corresponding sensitivity equations.
Cocquet et al.~\cite{Cocquet2021} approximate the unsteady SWE with a Helmholtz equation and use a subspace trust-region method for the minimisation. 
They consider a continuous optimisation problem for the weak form of their steady-state equation, prove continuity of the control-to-state operator, existence of an optimal solution and convergence ot the Finite Element approximation.
Sch{\"a}fer et al.~\cite{SchaeferEtAl2021} prove convergence of discretisation schemes for the adjoint equation of a hyperbolic conservation law with boundary control.
Riffo~\cite{Riffo2019} performs reconstruction using PDE-constrained optimisation but under assumptions for which the SWE reduce to a steady state Helmholtz problem.  
Losch and Wunsch~\cite{Losch2003} use a linearized shallow water model to reconstruct a bathymetry with a sensitivity approach.
In the dissertation of Wu~\cite{Wu2022}, coastal bathymetries are inferred from surface waves using a high-order spectral method for simulations and a variational data assimilation method for the reconstruction. In Chapter 3 of the dissertation a detailed list of references on bathymetry reconstruction is provided.

There are also a variety of approaches that do not use an adjoint problem.
Vasan et al.~\cite{VasanEtAl2021} consider two inverse problems: Bathymetry reconstruction from free-surface deviation and velocity as well as velocity deduction from surface deviation and approximate bathymetry. 
They use one reconstruction algorithm for both problems at once.
Hajduk et al.~\cite{Hajduk2020} discretise the SWE with (dis)continuous Galerkin approaches and use a pseudo-time stepping scheme to reconstruct the bathymetry.
Gessese and Sellier~\cite{Gessese2013,Gessese2012} use a one-shot technique, where they replace water depth over ground in the SWE by the surface elevation and solve numerically for the difference. 
Dr{\"a}hne et al.~\cite{Draehne2016} simulate long wave run-up on a beach with a DG method and validate with experimental data but do not consider the inverse problem.
Hao et al.~\cite{HaoEtAl2024} use machine learning based techniques developed for image inpainting to fill in bathymetry data from measurements at the coast or beach. 
Their approach can provide normalised root means square errors (NRMSE) of around 20\%. 
They also provide a detailed discussion of the relevance of bathymetry reconstruction. 
Kar and Guha~\cite{KarGuha2018} compute the bathymetry by taking time averages of velocities and the free surface elevation using the fact that a hill in the bathymetry creates a dip in the free surface in the case of sub-critical flow.
Xie et al.~\cite{XieEtAl2023} use a neural network to reconstruct the bathymetry from side scan measurements.
Liu et al.~\cite{LiuEtAl2022} use a deep-learning based surrogate model for shallow water. 
They solve the inverse problem with a gradient descent method and automatic differentiation instead of a continuous adjoint.

While there are different techniques and theoretical analyses, so far there seem to be only few studies of how the PDE-constrained optimisation approaches fare for real data, in particular with a continuous adjoint.
One notable exception is the work by Lacasta et al.~\cite{LacastaEtAl2018,LacastaEtAl2019}, who use a continuous adjoint of the one-dimensional SWE to reconstruct inlet boundary conditions from experimental data in a similar setting as ours, although their wave tank is significantly larger.
However, they treat the bathymetry as known and do not attempt to reconstruct it.

As Khan et al.~\cite{Khan2021} state in their 2021 paper on bathymetry reconstruction, ``while the theoretical results from such models are promising, their applicability to real-world measurements has not been established''.
To our knowledge, this paper is the first attempt at demonstrating that PDE-constrained optimisation with a continuous adjoint can be used to reconstruct a bathymetry from real-world measurements in an experimental setting.
Both the used code~\cite{angel_2024} and data~\cite{dataset} are available and allow for the reproduction of all results reported in the paper.

\section{Problem Description}
\label{sec:problem}
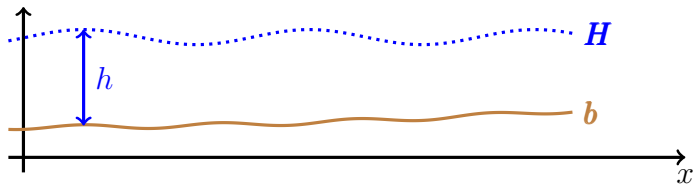
\begin{figure}[b]
	\centering
	\begin{tikzpicture}[domain=0:8, scale=1]
    x-axis
    \draw[->, very thick] (0,0) -- (9,0) node [below] {$x$};
    y-axis
    \draw[->, very thick] (.2,-.2) -- (.2,2);
    \draw[smooth, blue, dotted, very thick, samples=100, domain=0:7.5]
        plot(\x, {.1*cos(120*(\x-1))+1.6}) node [right] {\textcolor{blue}{$\pmb{H}$}};
    \draw[smooth, brown, very thick, samples=100, domain=0:7.5]
        plot(\x, {.03*cos(200*(\x-1))+.004*\x^2+0.4}) node [right] {$\pmb{b}$};
    \draw[<->, very thick, blue] (1,0.424) -- (1, 1.7) node[midway, right]{$h$};
	\end{tikzpicture}
	\caption{Free surface elevation $H$ (dotted blue), bathymetry $b$ (brown) and water height $h$. Our aim is to reconstruct $b$ from point measurements of $H$.}
	\label{fig:Hhb} 
\end{figure}

\begin{figure}[t]
	\centering
	\includegraphics[width=0.8\textwidth]{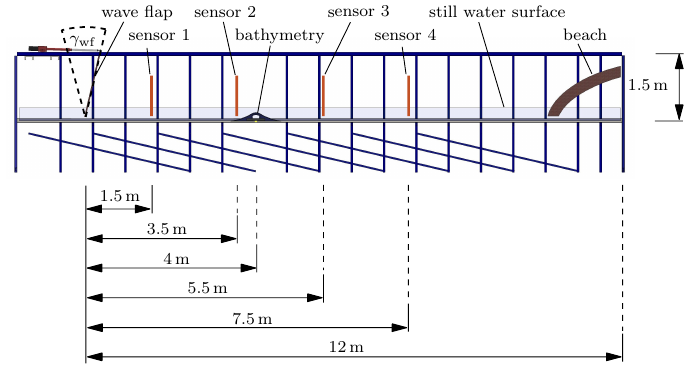}
	\caption{Sketch of used experimental setup using the wave flume at the Institute of Mechanics and Ocean Engineering at Hamburg University of Technology. Waves are generated at the left by a wave flap. Four sensors are positioned throughout the flume at increasing distance from the flap. The bathymetry is positioned between sensor 2 and sensor 3. At the right, a beach is installed to minimise reflections. In all simulations, we use the readings from sensor 1 to define the left boundary condition while readings from sensors 2, 3, 4 are used to reconstruct the bathymetry.}
	\label{fig:WaveFlumeSketch}
\end{figure}

Figure~\ref{fig:Hhb} shows a sketch of our problem setting.
The aim is to reconstruct the bathymetry $b$ from observations of the free surface elevation $H:= h+b$. 
We will study both simulated data and experimental measurements, generated by inserting a manufactured, roughly Gaussian-shaped bathymetry into a wave flume.
A sketch of the experimental setup is shown in Figure~\ref{fig:WaveFlumeSketch}.
A detailed explanation together with the data is available for download~\cite{dataset}.
The water height $h$ and the velocity $u$ of a homogeneous, incompressible fluid in the flume can be modelled by the nonlinear shallow water equations (SWE)
\begin{equation}\label{eq:nonconsvSWE}
 \begin{pmatrix}
  h\\
  u
 \end{pmatrix}_t
 +
 \begin{pmatrix}
  hu\\
  u^2 + gh
 \end{pmatrix}_x
 =
 \begin{pmatrix}
  0\\
  -gb_x - \kappa u
 \end{pmatrix}
\end{equation}
with a linear bottom friction term, where we set the corresponding coefficient $\kappa=0.2$. 
We experimented with more complex friction terms as well but found the impact on the accuracy of the model negligible.
The gravitational acceleration is set to $g=\SI{9.81}{m/s^2}$. 
Let $\Omega=[L, R]$ be the spatial domain and $Q:=\Omega \times [0,T]$ the space-time domain for some final time $T$. 

Suppose we have an observation $H_{\text{obs}}$ of the surface elevation. 
As we want to model waves in the flume, we need to carefully choose the boundary conditions such that the numerical solution of the SWE coincides with measurements obtained by the sensors. 
We use the observation $H_{\text{obs}}$ at sensor~1 to obtain the left boundary value for $h$ and thus set $L=\SI{1.5}{m}$, corresponding to the distance of sensor 1 from the wave flap, see Figure~\ref{fig:WaveFlumeSketch}.
As we use the data from the first sensor to set the left boundary condition, data from sensors 2 to 4 is left to be used in the objective function of the optimisation.

We observed good agreement between measurements and simulation when setting the right boundary condition for the velocity $u$ to zero and defining the right boundary to $R=\SI{15}{m}$ so that the computational domain is longer than the actual water flume.  
Setting the velocity at the boundary to form a simple wave can be beneficial for matching the experimental data when using finite volume methods~\cite{Draehne2016}, but realising this boundary condition in the used simulation framework proved to be unstable.
As the water is at rest at the start of the experiment, we set the initial velocity to zero.
In summary, the initial and boundary conditions read
\begin{subequations}
\begin{align}
	h(\cdot, 0) &= H_{\text{obs}}(\cdot, 0) - b \quad \text{on} \ \Omega \label{eq:swe_ic1}\\
	u(\cdot, 0) &= 0 \quad \text{on} \ \Omega \label{eq:swe_ic2}\\ 
	h(L, \cdot) &= H_{\text{obs}}(L, \cdot) - b(L) \quad \text{on} \ [0,T]\label{eq:swe_bc1} \\ 
	u(R, \cdot) &= 0 \quad \text{on} \  [0,T]. \label{eq:swe_bc2}
\end{align}
\end{subequations}
Note that in forward simulations during the reconstruction process, $b$ is the approximation to the exact bathymetry in the current iteration.

\subsection{Implementation}
The numerical solution of the SWE and its adjoint is computed with the open-source spectral solver framework Dedalus~\cite{Burns2020} written in Python. 
We used a Chebyshev basis as ansatz functions and a four-stage Runge Kutta IMEX method of third order included in the Dedalus package. 
To keep computation times reasonable, we chose a relatively coarse resolution with a grid size of $M=68$ and time step of $\SI{1e-3}{s}$.
The full code is available for download~\cite{angel_2024} and, together with the provided data~\cite{dataset}, can be used to reproduce all results reported in this paper.
\begin{figure}[ht]
	\centering
	\includegraphics[height=0.25\textheight]{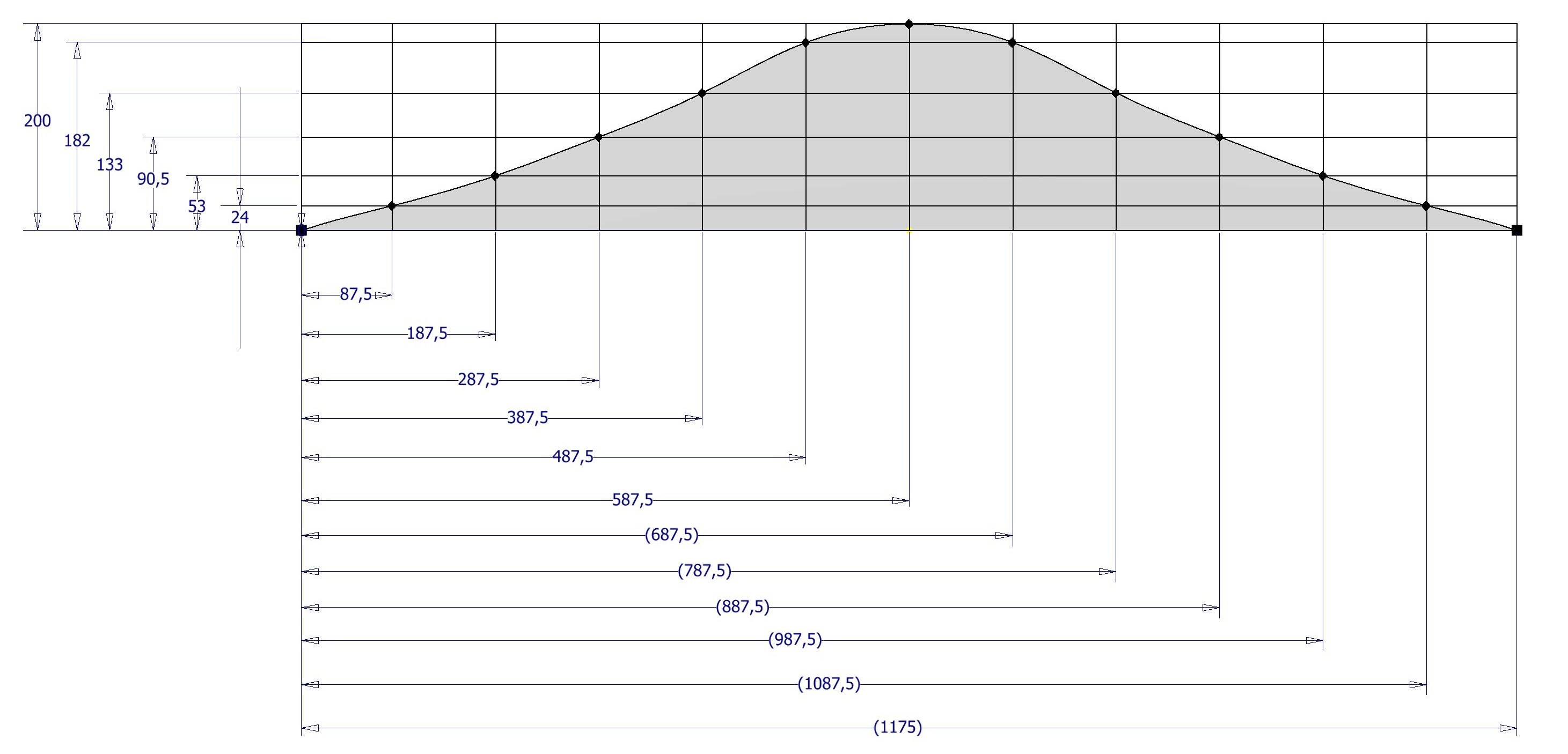}\\
	\hspace*{.8cm}\includegraphics[height=0.12\textheight]{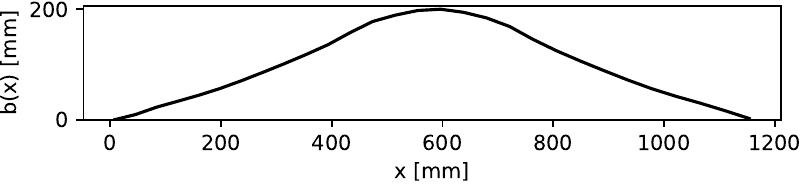}
	\caption{Measured points (top) and Dedalus representation (bottom) of the bathymetry that was used in the experiment. The Dedalus representation in the lower figure was generated by manually copying the coordinates from the data points in the upper sketch and connecting them using SciPy's \texttt{CubicSpline} function.}
	\label{fig:bathy}
\end{figure}

The representation of the bathymetry in Dedalus was generated by manually reading the data points in the technical sketch of the actual bathymetry in Fig.~\ref{fig:bathy} (top) and interpolating them using SciPy's \text{CubicSplines} function, resulting in the representation shown in  Fig.~\ref{fig:bathy} (bottom).
Note that the representation of the manufactured bathymetry in Figure~\ref{fig:bathy} was plotted with a much finer resolution than the one we used for the computations.
However, we emphasise that although the two bathymetries look similar, they are not identical.
This mismatch between the actual bathymetry and its representation in the code is part of the errors reported below, since the bathymetry reconstructed by our algorithm is compared against the representation in Dedalus and not the physical bathymetry.

\subsection{Validation of the forward model}
\begin{figure}[t]
	\centering
	\includegraphics[height=0.23\textheight]{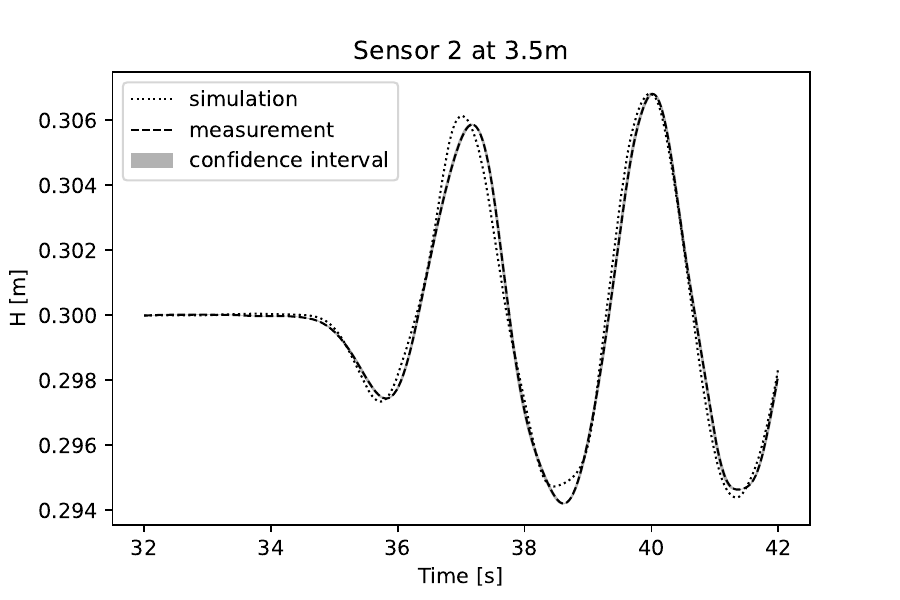}
	\includegraphics[height=0.23\textheight]{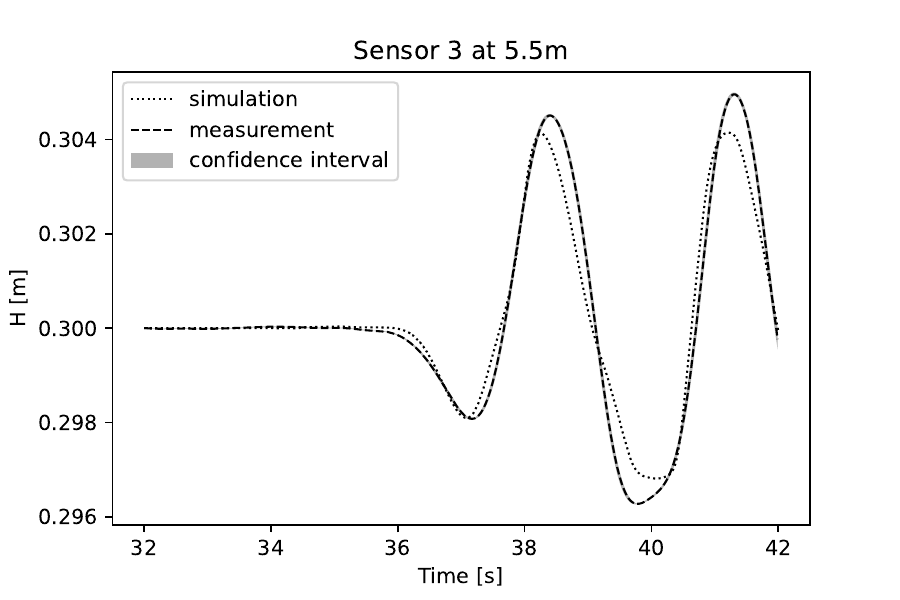}
	\includegraphics[height=0.23\textheight]{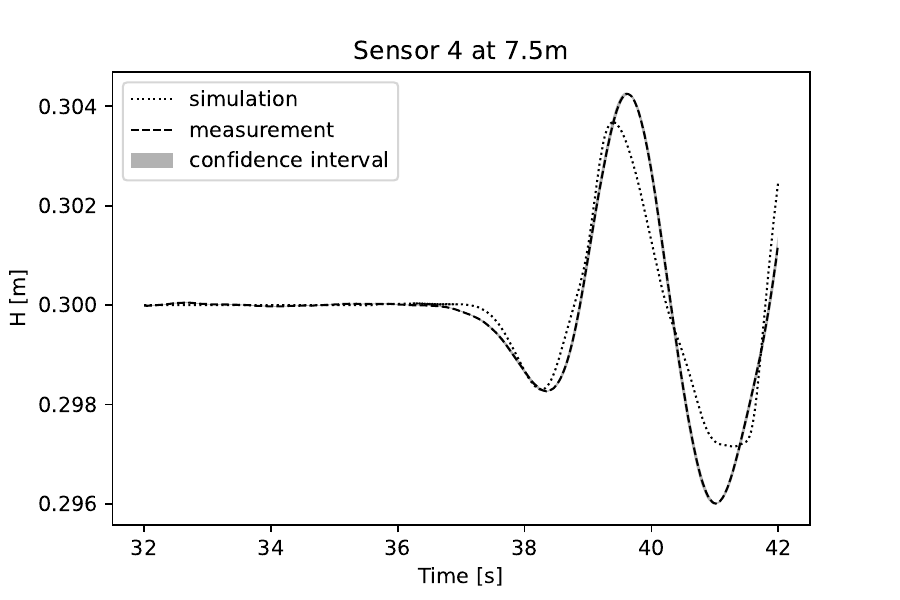}
	\caption{Measured (dashed) free surface elevation $H$ (dashed) with $95\%$ confidence interval (grey) and numerical simulation (dotted).}
	\label{fig:channel_forward_bathy}
\end{figure}

Figure~\ref{fig:channel_forward_bathy} shows the simulated surface elevation for the bathymetry sketched in Figure~\ref{fig:bathy} (bottom) as dotted line, the average of the $n=20$ measurements with bathymetry as dashed line and the two-sided 95 \% confidence interval~\cite{Dumbgen2016} as grey area for three sensors. For the computation of the confidence intervals we use the t-value $1.729$~\cite{tvalues} for $n-1$.
Note that because the generated waves change very little between the measurements, the confidence intervals are small and difficult to see in the plot.

Since measurements at sensor 1 provide the left boundary condition, we naturally have perfect agreement at $x=L$.
At sensor 2, model and experiment are in very good agreement with both phase and amplitude of the waves being very similar.
For sensors 3 and 4, we start to see an underestimation of amplitude in the model and a slightly too slow wave speed.
Most likely, this is due to dispersive effects not captured by the SWE.
While using a fully compressible model would likely generate a better match, it would also drastically increase computational cost.
Therefore, we stick to the relatively cheap SWE as model which, as we will show below, is already sufficient to provide a reasonable reconstruction of the bathymetry.
With $M=68$ and a time step of $\delta t= \SI{1e-3}{s}$, the relative $\ell^2$-errors between simulated and measured time series at sensor positions $2$, $3$, $4$ are around $1.7 \times 10^{-3}$, $1.7 \times 10^{-3}$ and $1.5 \times 10^{-3}$. 
Taking a finer spatial resolution with $M=100$ modes reduces these errors only marginally to $1.3 \times 10^{-3}$, $1.3 \times 10^{-3}$ and $1.1 \times 10^{-3}$.
We therefore use $M=68$ for all reconstructions for the sake of efficiency.

\begin{figure}[t]
	\centering
	\includegraphics[height=0.23\textheight]{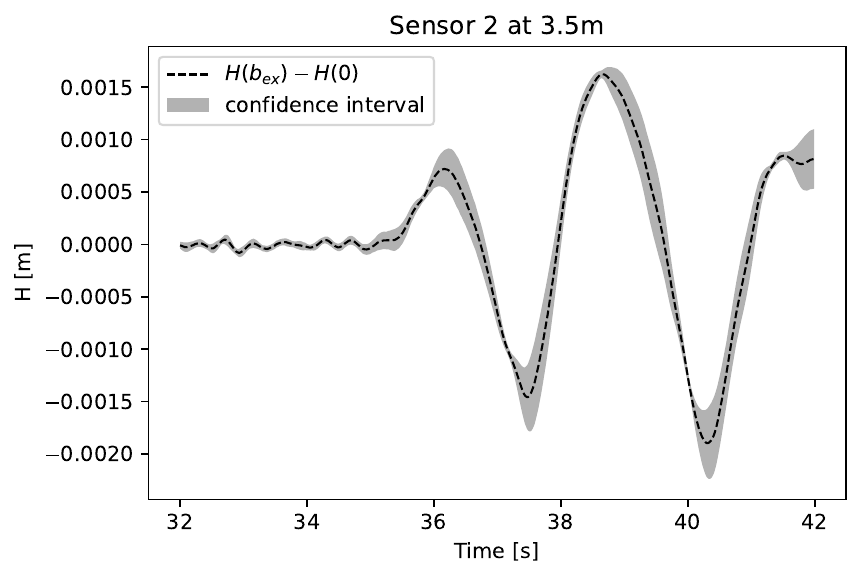}
	\includegraphics[height=0.23\textheight]{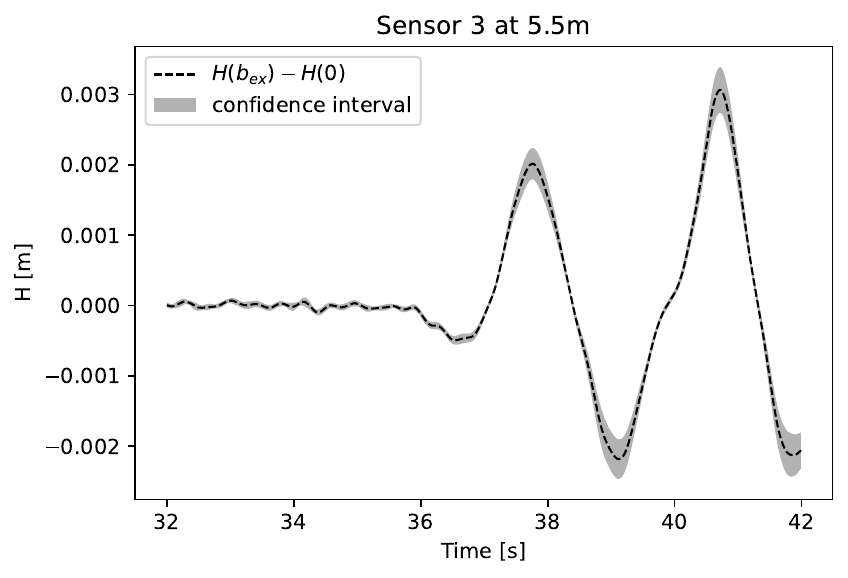}
	\includegraphics[height=0.23\textheight]{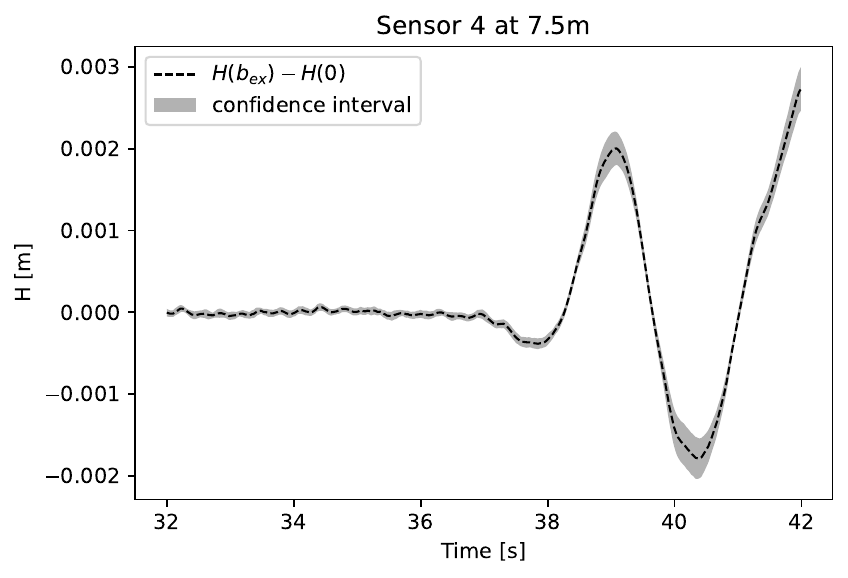}
	\caption{Difference between the averages of measurements with and without bathymetry (dashed) and $95\%$ confidence interval (grey).}
	\label{fig:diff_confidence}
\end{figure}

The information that can be used to reconstruct the bathymetry is the difference in the free surface elevation in the presence and absence of bathymetry. 
Figure~\ref{fig:diff_confidence} shows the difference between the average of the 20 experiments with and the average of the 20 experiments without bathymetry. 
Again, the $95\%$ confidence interval is shown in grey. 
Here we used the t-value $1.686$~\cite{tvalues} for $2(n-1)$, because we have a two-sample problem with two different experiments~\cite{Dumbgen2016}.
The maximum differences in the water height with and without bathymetry are around $\SI{1.5}{mm}$ at sensor 2 and $\SI{3}{mm}$ at sensors 3 and 4.

\section{Optimisation}
\label{sec:gradient_descent}
Bathymetry reconstruction can be formulated as a minimisation problem with the shallow water equations~\eqref{eq:nonconsvSWE} as constraint. 
We define a functional that measures the error between the simulated free surface elevation and the observation for a given bathymetry $b$. 
Let $\gamma, \delta, \lambda_1, \lambda_2 > 0$.
We define the objective functional
\begin{equation}\label{eq:J}
\begin{aligned}
	J(b,H):= &\frac{\gamma}{2} \int_Q
			(H-H_{\text{obs}})^2 ~d(x,t)
		+ \frac{\delta}{2} \int_{\Omega} (H(\cdot, T) -
			H_{\text{obs}}(\cdot, T))^2 ~dx \\
		&+ \frac{\lambda_1}{2}\int\limits_{\Omega} b^2 ~dx
		+ \frac{\lambda_2}{2} \int\limits_{\Omega} b_x^2 ~dx
\end{aligned}
\end{equation}
which consists of the $L^2$-error in the observation plus regularisation terms to ensure well-posedness of the PDE-constrained optimisation problem
\begin{equation}\label{eq:opt_problem}
 \min_{b \in H^1(\Omega)} J(b,H) \ \text{subject to} \ \eqref{eq:nonconsvSWE}.
\end{equation}
In this context, the SWE are called the forward model. 
To find a solution of~\eqref{eq:opt_problem}, we use gradient decent, and thus solve adjoint equations to compute the gradient. 
Since the gradient of the bathymetry $b_x$ appears in the SWE, we need to use the $H^1$-gradient. 
We obtain the adjoint equations using the Lagrange function~\cite{Troeltzsch2010} which in our case is defined as
\begin{equation}\label{eq:Lagrangian}
	\begin{aligned}
		\mathcal{L}(h,u,b,p) &= 
		J(b,H) - \int_Q (h_t + (hu)_x)p_1 ~d(x,t) \\
		&\quad - \int_Q \left( u_t + 2uu_x + g h_x + g b_x + \kappa u \right)p_2 ~d(x,t) \\
		&\quad - \int_0^T (h(L, \cdot) - H_{\text{obs}}(L, \cdot) + b(L))p_3 ~dt
		- \int_0^T u(R, \cdot)p_4 ~dt\\
		&\quad - \int_{\Omega} (h(\cdot, 0)-(H_{\text{obs}}(\cdot, 0)-b))p_5 ~dx
	\end{aligned}
\end{equation}
with Lagrange multipliers $p:=(p_1,\dots,p_5)$.
Finding the minimum of this function corresponds to solving~\eqref{eq:opt_problem} while partial derivatives of $\mathcal{L}$ yield the continuous adjoint equations and the $L^2$-gradient of the objective functional.

\subsection{Adjoint Problem and Gradient}
The adjoint equations for the SWE in quasi-linear form can be determined using the optimality condition
\begin{equation}\label{eq:opt_cond}
	\begin{aligned}
		\mathcal{L}_{h}(h,u,b,p)\varphi &= 0 \quad \forall \varphi \in C_0^{\infty}(Q)\\
		\mathcal{L}_{u}(h,u,b,p)\varphi &= 0 \quad \forall  \varphi \in C_0^{\infty}(Q).
	\end{aligned}
\end{equation}
The adjoint problem reads
\begin{equation}\label{eq:adjointProblem}
	\begin{aligned}
		&p_{1,t} + u p_{1,x} + g p_{2,x} = - \gamma \left( h + b - H_{\text{obs}} \right) \\
		&p_{2,t} + h p_{1,x} + 2up_{2,x} = \kappa p_2
	\end{aligned}
\end{equation}
with final and boundary conditions
\begin{subequations}
	\begin{align}
		&p_1(\cdot, T) = \delta \left( h(\cdot, 0) + b -
		H_{\text{obs}}(\cdot, 0) \right) \label{eq:adj_fc1}\\
		&p_2(\cdot, T) = 0 \label{eq:adj_fc2}\\
		&p_1(L, \cdot) = - 2\frac{u(L, \cdot)}{h(L, \cdot)} p_2(L, \cdot) \label{eq:adj_bc1}\\
		&p_2(R, \cdot) = 0, \label{eq:adj_bc2}
	\end{align}
\end{subequations}
see~\ref{app:derivation} for the detailed derivation.
Note that the division in~\eqref{eq:adj_bc1} is meant pointwise and allowed as $h(L,t)>0\ \forall t \in [0,T]$.
Let
\begin{equation}
 \tilde{v} := \int_0^T \gamma(h - H_{\text{obs}}) + g p_{2,x} ~dt + bT 
    + \delta \left( h(\cdot, T) + b - H_{\text{obs}}(\cdot, T) \right)
    + \lambda_1 b - \lambda_2 b_{xx} - p_5.
\end{equation}
Then the $H^1$-gradient is the solution $v$ of the problem
\begin{equation}\label{eq:poisson-like}
\begin{aligned}
 &v - \Delta v = \tilde{v} \\
 &v(R) = 0 \\
 &v(L) =0.
\end{aligned}
\end{equation}
By choosing homogeneous Dirichlet boundary conditions for~\eqref{eq:poisson-like} we assume that the bathymetry does not change at the boundaries. As we initialise the bathymetry with zero throughout all reconstruction experiments in this paper we have $b(L)=0$. Therefore, in~\eqref{eq:swe_bc1} we always have $h(L, \cdot) = H_{\text{obs}}(L, \cdot)$.
A different boundary condition for the gradient, however, could also model non-zero bathymetry at the domain boundary.
\subsection{Minimisation of the objective function}
We use gradient descent for minimisation as it is easy to implement. 
Given an observation of the surface elevation $H_{\text{obs}}$, it starts with an initial guess for the bathymetry $b$. 
Using an adaptive step size $a_j$ found by Armijo line search~\cite{Armijo1966}, the algorithm descends in the direction of the antigradient.
In each iteration of gradient descent the forward equations~\eqref{eq:nonconsvSWE} and the adjoint problem~\eqref{eq:adjointProblem} need to be solved. 
The results and the current bathymetry $b$ are then used to compute the gradient $v$. 
Our algorithm stops if the $\ell^2$-norm of the gradient is smaller than a user-defined tolerance $\varepsilon$.

\section{Reconstruction from simulated observations}\label{sec:sim_obs}
In this section, we analyse how well the reconstruction works for data that have been generated by the forward model.
This allows us to control noise, eliminates the effect of model error and provides a ``best case'' to compare reconstruction from experimental data against.
We use the Dedalus representation of the exact bathymetry shown in Figure~\ref{fig:bathy} (bottom). 
To generate the observation, we compute the solution $(h, u)$ of~\eqref{eq:nonconsvSWE} to~\eqref{eq:swe_bc2} for the exact bathymetry $b_{\text{ex}}$ and determine $H_{\text{obs}}=h+b_{\text{ex}}$ on the whole spatial domain. 
We used $M=100$ spatial grid points and a time step of $\SI{5e-5}{s}$ to compute $H_{\text{obs}}$, which is different from the resolution used in the reconstruction. 
Using exactly the same discretisation would ignore the discretisation error in the optimisation and is called an ``inverse crime'' that should be avoided~\cite{Kaipio2007}.
In order to deal with the different resolutions, Dedalus's \texttt{change\_scale} function is used to interpolate between spatial meshes and SciPy's \texttt{CubicSpline} function to interpolate in time.

Note that the bathymetries in this section are only shown for $x<\SI{10}{m}$ for better visibility. The oscillations in the interval $[\SI{10}{m}, \SI{15}{m}]$ are fairly small and thus for the plots negligible.

\subsection{Observation on the whole spatial domain}
The minimisation is performed as described in Section~\ref{sec:gradient_descent}, with a computed observation that is available on the whole spatial domain. We set the coefficients in the objective functional~\eqref{eq:J} to $\gamma=\delta=0.5$, $\lambda_1=10^{-6}$ and $\lambda_2=10^{-7}$.
For a second reconstruction from noisy data we added normally distributed noise to the computed observation with zero mean and a standard deviation of $5 \%$ of the highest observed local wave amplitude.
For this case, we increased the regularisation parameter to $\lambda_1=10^{-5}$.
\begin{figure}[b]
	\centering
	\includegraphics[scale=1]{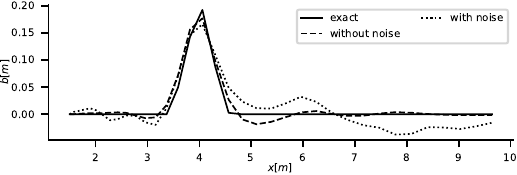}
	\includegraphics[scale=1]{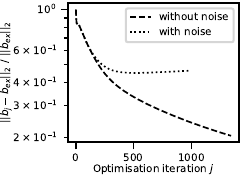}
	\caption{Reconstructed bathymetries from simulated data with and without noise (left) and relative $\ell^2$-error to exact bathymetry (right).}
	\label{fig:sim_everywhere}
\end{figure}
The corresponding reconstructions and relative errors to the exact bathymetry are shown in Figure~\ref{fig:sim_everywhere}.
Noise in the observation leads to some additional oscillations in the reconstructed bathymetry. In this case, the algorithm stopped early as there was no step size found for which the update of the bathymetry leads to a smaller value of the objective functional, indicating that the approximated anti-gradient was no longer a descent direction.

\subsection{Observation at sensor positions}
\begin{figure}[t]
	\centering
	\includegraphics[scale=1]{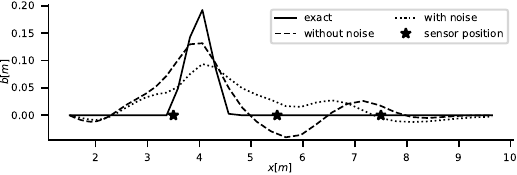}
	\includegraphics[scale=1]{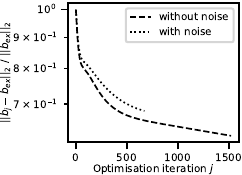}
	\caption{Reconstructed bathymetries from simulated data at sensor positions (left) and relative $\ell^2$-error to exact bathymetry (right).}
	\label{fig:sim_data234}
\end{figure}

In order to get closer to the experimental setup, we additionally reconstructed the bathymetry from observations only at three points in space where there are sensors in the real experiment.
Let $x_i$, $i=1,\dots,m_p$, be the sensor positions and $H_{\text{obs}, i}(t)$ the interpolated corresponding observations. 
To avoid using Dirac distributions in the objective functional, we replace the point measurements by Gaussians centred at the sensor positions, with maximum value equal to $H_{\text{obs}, i}$. 
The variance of the Gaussians is set to $\sqrt{0.045}$, making it narrow but still somewhat wider than the distance between two adjacent points on the grid.
This changes the objective functional~\eqref{eq:J} to
\begin{equation}\label{eq:J_gauss}
	\begin{aligned}
		J(b,H):= &\frac{\gamma}{2} \int_Q \left(\sum_{i=1}^{m_p} 	
		(H(x_i, t)-H_{\text{obs}, i}(t)) 
		\exp \left(- \frac{1}{2}\frac{(x-x_i)^2}{0.045} \right) \right)^2 ~d(x,t) \\
		&+ \frac{\delta}{2} \int_{\Omega} \left(\sum_{i=1}^{m_p}
		(H(x_i, T) - H_{\text{obs}, i}(T)) 
		\exp \left(- \frac{1}{2}\frac{(x-x_i)^2}{0.045} \right) \right)^2 ~dx \\
		&+ \frac{\lambda_1}{2}\int\limits_{\Omega} b^2 ~dx
		+ \frac{\lambda_2}{2} \int\limits_{\Omega} b_x^2 ~dx.
	\end{aligned}
\end{equation}
Figure~\ref{fig:sim_data234} shows the reconstruction for this case without noise as well as with noise added to the observation.
The parameters are the same as for the previous case with observation on the whole spatial domain. 
Since the algorithms has a lot less information to work with, the reconstructed bathymetry is less accurate.
However, with and without added noise, the maximum of the hill is approximately at the correct position, although the height is underestimated and the hill is somewhat smeared out.
There are also more spurious oscillations between sensor 3 and sensor 4.

\section{Reconstruction from experimental observations}
\begin{figure}
    \centering
    \includegraphics[scale=1]{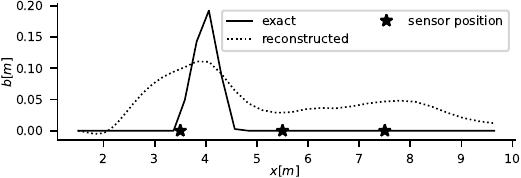}
    \includegraphics[scale=1]{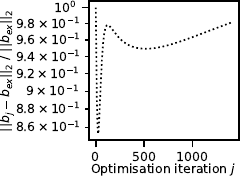}
    \caption{Reconstructed bathymetries from measurements.}
    \label{fig:reconstr_meas}
\end{figure}

For the reconstruction from experimental measurements we used the same parameters as for the optimisation from simulated observation at sensor positions with added noise in Section~\ref{sec:sim_obs}. 
Figure~\ref{fig:reconstr_meas} shows the reconstruction from the measurement data obtained by sensors 2, 3 and 4. 
The maximum of the reconstructed hill is at the correct position, but second hill at $x \approx \SI{8}{m}$ has appeared. 
Compared to the reconstruction from the simulated observation with added noise, the unwanted hill on the right is larger, but the maximum of the hill at $\SI{4}{m}$ is closer to the exact bathymetry. 
In order to measure the quality of the reconstructions, we compare the $\ell^2$- and $\ell^{\infty}$- errors as well as the normalised root mean square error (NRMSE)
\begin{equation}
	\text{NRMSE} = \frac{\sqrt{\frac{1}{M} \sum_{i=1}^{M} b_i-b_{\text{ex}, i}}}{b_{\text{ex, max}} - b_{\text{ex, min}}}
\end{equation}
with $M$ being the number of points where the representation of the exact bathymetry $b_{\text{ex}}$ is defined. In our case, we have $M=68$.

We show relative errors and NRMSEs of different configurations for simulated and experimental observations in Table~\ref{tab:errors}. 
For simulated observation, where model errors are absent, errors increase when using less data.
However, using only sensors 2 and 3 provides a reconstruction that is almost as accurate as using all three sensors, indicating that sensor 4 contributes little useful information.
Adding noise increases the $\ell^2$- and $\ell^{\infty}$-error but has little impact on the NRMSE.

For real measurements, $\ell^{\infty}$-errors are comparable to those for simulated observations while $\ell^{2}$-errors and NRMSE are higher.
In contrast to simulated observations, the lowest errors are achieved when using only sensors 2 and 3 for reconstruction.
This is likely because of the lack of dispersive effects in the SWE and the resulting mismatch between reality and model at sensor 4.
Interestingly, Lacasta et al.~\cite{LacastaEtAl2019} also report issues with the shallow water model, although in their case the reason is the formation of surf.
To put these numbers into context, for their state-of-the-art ML-based approach, Hao et al. report NRMSE between 20.6\% and 31.8\%~\cite[Table 2]{HaoEtAl2024} for real-world test cases.
Considering that our setup is simpler, the NRMSE of about 14\% for measured observations in Table~\ref{tab:errors} seems reasonable.
\begin{table}[ht]
	\centering
	\begin{tabular}{@{}llrrr@{}} \toprule
	\multicolumn{2}{c}{Observation} \\ \cmidrule(r){1-2}
	Type & Used sensors & $\ell^2$-error & $\ell^{\infty}$-error & NRMSE \\ \midrule
	Simulated  & 2,3,4 & 62.11 & 37.51 & 10.17 \\
	& 2,3 & 64.45 & 37.73 & 10.55\\
	& 2 & 70.17 & 46.63 & 11.48\\
	& 3 & 91.88 & 69.27 & 15.04 \\
	Simulated + noise & 2,3,4 & 68.20 & 51.34 & 11.16\\
	& 2,3 & 68.32 & 49.13 & 11.18\\
	& 2 & 66.02 & 52.14 & 10.81\\
	& 3 & 87.31 & 78.7 & 14.29 \\
	Measured & 2,3,4 & 97.97 & 49.06 & 16.07\\
	& 2,3 & 84.94 & 48.98 & 13.9\\
	& 2 & 92.76 & 53.38 & 15.18\\
	& 3 & 98.69 & 61.5 & 16.15 \\ \bottomrule
	\end{tabular}
	\caption{Relative error in the $\ell^2$- and $\ell^{\infty}$-norm and normalised root mean squared error (NRMSE) between reconstructed bathymetry and the model representation of the exact bathymetry for reconstructions using data from different sensors.}
	\label{tab:errors}
\end{table}

\section{Conclusions}
This paper analyses the performance of a PDE-constrained optimisation approach with a continuous adjoint for the reconstruction of bathymetry from experimental data.
Our setup consists of a roughly Gaussian shaped bathymetry placed in a wave flume with four sensors providing point measurements over time.
Data from the first sensor is used to prescribe the boundary condition for the model, leaving measurements from the remaining three sensors for the reconstruction.

To establish a baseline without the influences of modelling errors, we first perform the reconstruction using simulated observations with and without added noise.
Then, we perform the reconstruction from experimental data using data from all three sensors, two sensors or only a single sensor.
Even though the differences in wave heights with and without bathymetry are only a few millimetres, the approach can provide a qualitatively correct reconstruction of the bathymetry.
The maximum of the bathymetry is placed approximately at the right position, even though its height is underestimated and the reconstruction is somewhat smeared out compared to the original.
Achieved normalised root mean square errors (NRMSE) of around 14\% are found to be in line with recent results from the literature for a machine learning based approach.
Our results also suggest that the quality of the reconstruction is sensitive to sensor placement.
In theory, the optimisation approach could also be modified to also adjust sensor placement but investigating this is left for future work.

\section*{Acknowledgments}
This project has received funding from the European High-Performance Computing Joint Undertaking (JU) under grant agreement No 955701. The JU receives support from the European Union’s Horizon 2020 research and innovation programme and Belgium, France, Germany, and Switzerland. This project also received funding from the German Federal Ministry of Education and Research (BMBF) grant 16HPC048.
We acknowledge the support by the Deutsche Forschungsgemeinschaft (DFG) within the Research Training Group GRK 2583 ``Modeling, Simulation and Optimization of Fluid Dynamic Applications''. This project was also supported by the DFG via grant number 528383251.

Our special thanks go to Riza Demir and Wolfgang Brennecke of the Institute of Mechanics and Ocean Engineering at TUHH, who built and measured the bathymetry. 
We also thank Vanessa Trapp from the Chair for Stochastics at TU Hamburg for providing us with references regarding the computation of confidence intervals.

\appendix
\section{Derivation of Adjoint Equations and Gradient}\label{app:derivation}
Continuous adjoints to the nonlinear SWE have been derived before~\cite{Sanders2000,LacastaEtAl2018, LacastaEtAl2019}, but for the conservative form of the SWE and a different bottom friction term. \
Here we derive adjoint equations for the quasi-linear form. The Lagrangian for our problem reads
\begin{equation}
	\begin{aligned}
		\mathcal{L}(h,u,b,p) &= 
		\frac{\gamma}{2} \int_{Q} (h + b - H_{\text{obs}})^2 ~d(x,t) \\
		&\quad + \frac{\delta}{2} \int_{\Omega} (h(\cdot, T) + b -   
		H_{\text{obs}}(\cdot, T))^2 ~dx
		+ \frac{\lambda}{2} \int_{\Omega} b_x^2 ~dx \\
		&\quad - \int_Q (h_t + (hu)_x)p_1 ~d(x,t) \\
		&\quad - \int_Q \left( u_t + 2uu_x + g h_x + g b_x + \kappa u \right)p_2 ~d(x,t) \\
		&\quad - \int_0^T (h(L, \cdot) - H_{\text{obs}}(L, \cdot) + b(L))p_3 ~dt
		- \int_0^T u(R, \cdot)p_4 ~dt\\
		&\quad - \int_{\Omega} (h(\cdot, 0)-(H_{\text{obs}}(\cdot, 0)-b))p_5 ~dx.
	\end{aligned}
\end{equation}
with Lagrange multipliers $p_1, \dots, p_5$.
We use the necessary optimality condition
\begin{equation}\label{eq:opt1}
	\nabla_{(h,u)}\mathcal{L}(h,u,b,p) = 0
\end{equation}
to determine the continuous adjoint. Using integration by parts we get
\begin{equation}\label{eq:partial_h}
	\begin{aligned}
		\mathcal{L}_{h}(h,u,b,p)\varphi &= 
		\gamma \int_Q \left( h+b-H_{\text{obs}} \right)\varphi ~d(x,t) \\
		&\quad + \delta \int_{\Omega} \left( h(\cdot, T)+b-  
		H_{\text{obs}}(\cdot, T) \right)\varphi(\cdot,T) ~dx \\
		&\quad - \int_Q \varphi_t p_1 ~d(x,t)
		- \int_Q (\varphi u)_x p_1 ~d(x,t)
		- \int_Q g \varphi_x p_2 ~d(x,t) \\
		&\quad - \int_0^T \varphi(L, \cdot)p_3 ~dt
		- \int_{\Omega} \varphi(\cdot, 0) p_5 ~dx \\
		&= \gamma \int_Q \left( h+b-H_{\text{obs}} \right)\varphi ~d(x,t) \\
		&\quad + \delta \int_{\Omega} \left( h(\cdot,T)+b-H_{\text{obs}}(\cdot,T) \right)\varphi(\cdot,T) ~dx \\
		&\qquad + \int_Q \varphi p_{1,t} ~d(x,t) - \int_{\Omega} \left[ \varphi p_1 \right]_0^T ~dx
		+ \int_Q \varphi u p_{1,x} ~d(x,t) \\
		&\qquad - \int_0^T \left[\varphi u p_1 \right]_L^R ~dt
		+ \int_Q g \varphi p_{2,x} ~d(x,t) 
		- \int_0^T \left[g \varphi p_2 \right]_L^R ~dt  \\
		&\qquad   - \int_0^T \varphi(L, \cdot) p_3 ~dt
		- \int_{\Omega} \varphi(\cdot, 0) p_5 ~dx,
	\end{aligned}
\end{equation}
which has to be zero for all $\varphi \in C_0^{\infty}(Q)$ due to the optimality condition~\eqref{eq:opt_cond}. We obtain
\begin{equation*}
	\int_Q \left( \gamma \left( h+b-H_{\text{obs}} \right) 
	+ p_{1,t} + u p_{1,x} + g p_{2,x} \right)\varphi ~d(x,t) = 0
\end{equation*}
and, using the fundamental lemma of calculus of variations, we get the first part of the adjoint problem
\begin{equation}\label{eq:nonconsv_adj1}
	p_{1,t} + u p_{1,x} + g p_{2,x} = - \gamma \left( h + b - H_{\text{obs}} \right).
\end{equation}
Now we let single conditions to the function $\varphi$ vary in order to determine terminal and boundary conditions for the adjoint equation.
By allowing any value for $\varphi(\cdot, T)$ in Equation~\eqref{eq:partial_h} we get the final time condition for the first adjoint variable
\begin{equation}
	p_1(\cdot, T) = \delta \left( h(\cdot, T) + b - H_{\text{obs}}(\cdot, T) \right).
\end{equation}
Without the condition $\varphi(\cdot, 0)=0$ we obtain
\begin{equation}
	p_1(\cdot, 0) = p_5
\end{equation}
and letting $\varphi(R, \cdot)$ vary we get
\begin{equation}\label{eq:p2R0}
	g p_2(R, \cdot) = - u(R, \cdot) p_1(R, \cdot) = 0.
\end{equation}
Dropping the condition $\varphi(L, \cdot)=0$ yields
\begin{equation}\label{eq:primary_p3}
	p_3 = - u(L, \cdot) p_1(L, \cdot) - g p_2(L, \cdot)
\end{equation}
for the Lagrange multiplier $p_3$.
The next partial derivative is
\begin{equation}\label{eq:partial_u}
	\begin{aligned}
		\mathcal{L}_{u}(h,u,b,p)\varphi &= - \int_Q (h \varphi)_x p_1 ~d(x,t)
		- \int_Q \left( \varphi_t
		+ 2 \left( u \varphi \right)_x + \kappa \varphi \right)p_2 ~d(x,t) \\
		&= \int_Q h \varphi p_{1,x} ~d(x,t)
		- \int_0^T \left[ h \varphi p_1 \right]_L^R ~dt
		+ \int_Q \varphi p_{2,t} ~d(x,t) \\
		&\quad - \int_{\Omega} \left[ \varphi p_2 \right]_0^T ~dx
		+ \int_Q 2 u \varphi p_{2,x} ~d(x,t) \\
		&\quad - \int_0^T \left[ 2u \varphi p_2 \right]_L^R ~dt
		- \int_Q \kappa \varphi p_2 ~d(x,t)
		- \int_0^T\varphi(R, \cdot)p_4 ~dt .
	\end{aligned}
\end{equation}
Hence, the second part of the adjoint equation is
\begin{equation}\label{eq:nonconsv_adj2}
	p_{2,t} + h p_{1,x} + 2up_{2,x} = \kappa p_2.
\end{equation}
We let $\varphi(\cdot, T)$ vary and obtain the terminal condition for the second adjoint variable
\begin{equation}
	p_2(\cdot, T) = 0.
\end{equation}
If we leave out the condition $\varphi(L, \cdot)=0$, we get
\begin{equation}
	\begin{aligned}
		2u(L, \cdot) p_2(L, \cdot) &= - p_1(L, \cdot) h(L, \cdot) \\
		\Leftrightarrow p_1(L, \cdot) &= - 2\frac{u(L, \cdot)}{h(L, \cdot)} p_2(L, \cdot).
	\end{aligned}
\end{equation}
Together with~\eqref{eq:primary_p3} this implies that
\begin{equation}\label{eq:p3}
 p_3 = \left( 2 \frac{(u(L, \cdot))^2}{h(L, \cdot)} - g \right) p_2(L, \cdot).
\end{equation}
To rewrite the adjoint problem as an initial value problem we set $\tau:=T-t$ and define
\begin{align*}
	\tilde{p}(x, \tau) &:= p(x, T-t) \\
	\tilde{h}(x, \tau) &:= h(x, T-t) \\
	\tilde{u}(x, \tau) &:= u(x, T-t).
\end{align*}
The time transformed adjoint problem reads
\begin{equation}
	\begin{aligned}
		&\tilde{p}_{1,\tau} - \tilde{u} \tilde{p}_{1,x} - g \tilde{p}_{2,x} = \gamma \left( \tilde{h} + b - \tilde{H}_{\text{obs}} \right) \\
		&\tilde{p}_{2,\tau} - \tilde{h} \tilde{p}_{1,x} - 2\tilde{u}\tilde{p}_{2,x} = - \kappa \tilde{p}_2 \\
		&\tilde{p}_1(\cdot, 0) =
		 \delta \left( \tilde{h}(\cdot, 0) + b -
		  \tilde{H}_{\text{obs}}(\cdot, 0) \right) \\
		&\tilde{p}_2(\cdot, 0) = 0 \\
		&\tilde{p}_1(L, \cdot) =
		 - 2\frac{\tilde{u}(L, \cdot)}{\tilde{h}(L, \cdot)} \tilde{p}_2(L, \cdot) \\
		&\tilde{p}_2(R, \cdot) = 0.
	\end{aligned}
\end{equation}
The $L^2$-gradient is
\begin{equation}\label{eq:partial_b}
	\begin{aligned}
		\mathcal{L}_{b}(h,u,b,p)\delta b &= \int_{\Omega} \left( 
		\int_0^T \gamma(h - H_{\text{obs}}) + g p_{2,x} ~dt + bT \right) \delta b ~dx \\
		&+ \int_{\Omega} \left( \delta \left( h(\cdot, T) + b - H_{\text{obs}}(\cdot, T) \right)
		- \lambda b_{xx}
		- p_5 \right) \delta b ~dx =: \tilde{v}.
	\end{aligned}
\end{equation}

\bibliographystyle{elsarticle-num} 
\bibliography{refs}

\begin{thebibliography}{10}
\expandafter\ifx\csname url\endcsname\relax
  \def\url#1{\texttt{#1}}\fi
\expandafter\ifx\csname urlprefix\endcsname\relax\def\urlprefix{URL }\fi
\expandafter\ifx\csname href\endcsname\relax
  \def\href#1#2{#2} \def\path#1{#1}\fi

\bibitem{Rabinovich2020}
A.~Rabinovich, \href{https://doi.org/10.1007/s00024-019-02349-3}{Twenty-seven
  years of progress in the science of meteorological tsunamis following the
  1992 daytona beach event}, Pure and Applied Geophysics 177 (2020) 1193--1230.
\newblock \href {https://doi.org/10.1007/s00024-019-02349-3}
  {\path{doi:10.1007/s00024-019-02349-3}}.
\newline\urlprefix\url{https://doi.org/10.1007/s00024-019-02349-3}

\bibitem{Sanders2000}
B.~Sanders, N.~Katopodes, Adjoint sensitivity analysis for shallow-water wave
  control, Journal of Engineering Mechanics-asce - J ENG MECH-ASCE 126 (09
  2000).
\newblock \href {https://doi.org/10.1061/(ASCE)0733-9399(2000)126:9(909)}
  {\path{doi:10.1061/(ASCE)0733-9399(2000)126:9(909)}}.

\bibitem{smith2004conventional}
W.~H.~F. Smith, D.~T. Sandwell,
  \href{https://doi.org/10.5670/oceanog.2004.63}{Conventional bathymetry,
  bathymetry from space, and geodetic altimetry}, Oceanography 17~(1) (2004)
  8--23.
\newblock \href {https://doi.org/10.5670/oceanog.2004.63}
  {\path{doi:10.5670/oceanog.2004.63}}.
\newline\urlprefix\url{https://doi.org/10.5670/oceanog.2004.63}

\bibitem{sagawa2019satellite}
T.~Sagawa, Y.~Yamashita, T.~Okumura, T.~Yamanokuchi,
  \href{https://doi.org/10.3390/rs11101155}{Satellite derived bathymetry using
  machine learning and multi-temporal satellite images}, Remote Sensing 11~(10)
  (2019) 1155.
\newblock \href {https://doi.org/10.3390/rs11101155}
  {\path{doi:10.3390/rs11101155}}.
\newline\urlprefix\url{https://doi.org/10.3390/rs11101155}

\bibitem{irish1999scanning}
J.~L. Irish, W.~J. Lillycrop,
  \href{https://www.sciencedirect.com/science/article/pii/S0924271699000039}{Scanning
  laser mapping of the coastal zone: the {SHOALS} system}, ISPRS Journal of
  Photogrammetry and Remote Sensing 54~(2-3) (1999) 123--129.
\newblock \href {https://doi.org/10.1016/S0924-2716(99)00003-9}
  {\path{doi:10.1016/S0924-2716(99)00003-9}}.
\newline\urlprefix\url{https://www.sciencedirect.com/science/article/pii/S0924271699000039}

\bibitem{ashphaq2021review}
M.~Ashphaq, P.~K. Srivastava, D.~Mitra,
  \href{https://doi.org/10.1016/j.joes.2021.02.006}{Review of near-shore
  satellite derived bathymetry: {C}lassification and account of five decades of
  coastal bathymetry research}, Journal of Ocean Engineering and Science 6~(4)
  (2021) 340--359.
\newblock \href {https://doi.org/10.1016/j.joes.2021.02.006}
  {\path{doi:10.1016/j.joes.2021.02.006}}.
\newline\urlprefix\url{https://doi.org/10.1016/j.joes.2021.02.006}

\bibitem{evagorou2022evaluation}
E.~Evagorou, A.~Argyriou, N.~Papadopoulos, C.~Mettas, G.~Alexandrakis,
  D.~Hadjimitsis, \href{https://www.mdpi.com/2072-4292/14/3/772}{Evaluation of
  satellite-derived bathymetry from high and medium-resolution sensors using
  empirical methods}, Remote Sensing 14~(3) (2022) 772.
\newblock \href {https://doi.org/10.3390/rs14030772}
  {\path{doi:10.3390/rs14030772}}.
\newline\urlprefix\url{https://www.mdpi.com/2072-4292/14/3/772}

\bibitem{BeckersEtAl2019}
S.~Beckers, J.~Behrens, W.~Wollner,
  \href{https://dx.doi.org/10.1016/j.apnum.2019.05.016}{Duality based error
  estimation in the presence of discontinuities}, Applied Numerical Mathematics
  144 (2019) 83–99.
\newblock \href {https://doi.org/10.1016/j.apnum.2019.05.016}
  {\path{doi:10.1016/j.apnum.2019.05.016}}.
\newline\urlprefix\url{https://dx.doi.org/10.1016/j.apnum.2019.05.016}

\bibitem{Khan2021}
R.~A. Khan, N.~K.-R. Kevlahan, {Variational assimilation of surface wave data
  for bathymetry reconstruction. Part I: algorithm and test cases}, Tellus A:
  Dynamic Meteorology and Oceanography (Jan 2021).
\newblock \href {https://doi.org/10.1080/16000870.2021.1976907}
  {\path{doi:10.1080/16000870.2021.1976907}}.

\bibitem{Khan2022}
R.~A. Khan, N.~K.-R. Kevlahan, {Variational Assimilation of Surface Wave Data
  for Bathymetry Reconstruction. Part II: Second Order Adjoint Sensitivity
  Analysis}, Tellus A: Dynamic Meteorology and Oceanography (Apr 2022).
\newblock \href {https://doi.org/10.16993/tellusa.36}
  {\path{doi:10.16993/tellusa.36}}.

\bibitem{Cocquet2021}
P.-H. Cocquet, S.~Riffo, J.~Salomon,
  \href{https://doi.org/10.1137/20M1326337}{{Optimization of Bathymetry for
  Long Waves with Small Amplitude.}}, SIAM J. Control. Optim. 59~(6) (2021)
  4429--4456.
\newblock \href {https://doi.org/10.1137/20M1326337}
  {\path{doi:10.1137/20M1326337}}.
\newline\urlprefix\url{https://doi.org/10.1137/20M1326337}

\bibitem{SchaeferEtAl2021}
P.~S. Aguilar, S.~Ulbrich,
  \href{https://opus4.kobv.de/opus4-trr154/frontdoor/index/index/year/2021/docId/481}{Convergence
  of numerical adjoint schemes arising from optimal boundary control problems
  of hyperbolic conservation laws} (2021).
\newline\urlprefix\url{https://opus4.kobv.de/opus4-trr154/frontdoor/index/index/year/2021/docId/481}

\bibitem{Riffo2019}
S.~R. Riffo, Mathematical methods for marine energy extraction, Ph.D. thesis,
  École Doctorale de Dauphine (2021).

\bibitem{Losch2003}
M.~Losch, C.~Wunsch, Bottom topography as a control variable in an ocean model,
  Journal of atmospheric and oceanic technology, Vol. 20, Nr. 11, pages (2003)
  1685--1696.

\bibitem{Wu2022}
J.~Wu, Adjoint-based data assimilation for ocean waves, Ph.D. thesis,
  University of Minnesota (2022).

\bibitem{VasanEtAl2021}
V.~Vasan, Manisha, D.~Auroux,
  \href{https://doi.org/10.1111/sapm.12418}{Ocean-depth measurement using
  shallow-water wave models}, Studies in Applied Mathematics 147~(4) (2021)
  1481--1518.
\newblock \href {https://doi.org/10.1111/sapm.12418}
  {\path{doi:10.1111/sapm.12418}}.
\newline\urlprefix\url{https://doi.org/10.1111/sapm.12418}

\bibitem{Hajduk2020}
H.~Hajduk, D.~Kuzmin, V.~Aizinger, {Bathymetry Reconstruction Using Inverse
  ShallowWater Models: Finite Element Discretization and Regularization}, in:
  {Lecture Notes in Computational Science and Engineering}, Springer
  International Publishing, 2020, pp. 223--230.
\newblock \href {https://doi.org/10.1007/978-3-030-30705-9_20}
  {\path{doi:10.1007/978-3-030-30705-9_20}}.

\bibitem{Gessese2013}
A.~F. Gessese, Algorithms for bed topography reconstruction in geophysical
  flows, Ph.D. thesis, University of Canterbury (2013).

\bibitem{Gessese2012}
A.~Gessese, M.~Sellier, {A Direct Solution Approach to the Inverse
  Shallow-Water Problem}, Mathematical Problems in Engineering 2012 (2012)
  1--18.
\newblock \href {https://doi.org/10.1155/2012/417950}
  {\path{doi:10.1155/2012/417950}}.

\bibitem{Draehne2016}
U.~Dr{\"a}hne, N.~Goseberg, S.~Vater, N.~Beisiegel, J.~Behrens,
  \href{https://www.mdpi.com/2077-1312/4/1/1}{{An Experimental and Numerical
  Study of Long Wave Run-Up on a Plane Beach}}, Journal of Marine Science and
  Engineering 4~(1) (2016).
\newblock \href {https://doi.org/10.3390/jmse4010001}
  {\path{doi:10.3390/jmse4010001}}.
\newline\urlprefix\url{https://www.mdpi.com/2077-1312/4/1/1}

\bibitem{HaoEtAl2024}
Z.~Hao, F.~Chen, X.~Jia, X.~Cai, C.~Yang, Y.~Du, F.~Ling,
  \href{https://doi.org/10.1029/2023WR035781}{Grdl: A new global reservoir
  area-storage-depth data set derived through deep learning-based bathymetry
  reconstruction}, Water Resources Research 60~(1) (2024) e2023WR035781.
\newblock \href {https://doi.org/10.1029/2023WR035781}
  {\path{doi:10.1029/2023WR035781}}.
\newline\urlprefix\url{https://doi.org/10.1029/2023WR035781}

\bibitem{KarGuha2018}
S.~Kar, A.~Guha, \href{https://doi.org/10.1063/1.5055944}{Letter: Ocean
  bathymetry reconstruction from surface data using hydraulics theory}, Physics
  of Fluids 30~(12) (2018) 121701.
\newblock \href {https://doi.org/10.1063/1.5055944}
  {\path{doi:10.1063/1.5055944}}.
\newline\urlprefix\url{https://doi.org/10.1063/1.5055944}

\bibitem{XieEtAl2023}
Y.~Xie, N.~Bore, J.~Folkesson, Neural network normal estimation and bathymetry
  reconstruction from sidescan sonar, IEEE Journal of Oceanic Engineering
  48~(1) (2023) 218--232.
\newblock \href {https://doi.org/10.1109/JOE.2022.3194899}
  {\path{doi:10.1109/JOE.2022.3194899}}.

\bibitem{LiuEtAl2022}
X.~Liu, Y.~Song, C.~Shen, \href{https://arxiv.org/abs/2203.02821}{Bathymetry
  inversion using a deep-learning-based surrogate for shallow water equations
  solvers} (2022).
\newblock \href {http://arxiv.org/abs/2203.02821} {\path{arXiv:2203.02821}}.
\newline\urlprefix\url{https://arxiv.org/abs/2203.02821}

\bibitem{LacastaEtAl2018}
P.~B. {Asier Lacasta}, Mario Morales-Hernández, P.~García-Navarro,
  \href{https://doi.org/10.1080/00221686.2017.1300196}{Application of an
  adjoint-based optimization procedure for the optimal control of internal
  boundary conditions in the shallow water equations}, Journal of Hydraulic
  Research 56~(1) (2018) 111–123.
\newblock \href {https://doi.org/10.1080/00221686.2017.1300196}
  {\path{doi:10.1080/00221686.2017.1300196}}.
\newline\urlprefix\url{https://doi.org/10.1080/00221686.2017.1300196}

\bibitem{LacastaEtAl2019}
A.~Lacasta, D.~Caviedes-Voullieme, P.~García-Navarro,
  \href{https://doi.org/10.1007/978-3-319-89988-6_10}{Application of the
  Adjoint Method for the Reconstruction of the Boundary Condition in Unsteady
  Shallow Water Flow Simulation}, Springer International Publishing, 2019, p.
  157–172.
\newblock \href {https://doi.org/10.1007/978-3-319-89988-6_10}
  {\path{doi:10.1007/978-3-319-89988-6_10}}.
\newline\urlprefix\url{https://doi.org/10.1007/978-3-319-89988-6_10}

\bibitem{angel_2024}
J.~Angel, D.~Ruprecht, S.~Götschel,
  \href{https://doi.org/10.5281/zenodo.10848185}{Waterchannel} (Mar. 2024).
\newblock \href {https://doi.org/10.5281/zenodo.10848185}
  {\path{doi:10.5281/zenodo.10848185}}.
\newline\urlprefix\url{https://doi.org/10.5281/zenodo.10848185}

\bibitem{dataset}
J.~Angel, J.~Behrens, S.~Götschel, M.~Hollm, D.~Ruprecht, R.~Seifried,
  \href{https://doi.org/10.15480/882.9403}{Data artefact: Bathymetry
  reconstruction from experimental data with pde-constrained optimisation}
  (2024).
\newblock \href {https://doi.org/https://doi.org/10.15480/882.9403}
  {\path{doi:https://doi.org/10.15480/882.9403}}.
\newline\urlprefix\url{https://doi.org/10.15480/882.9403}

\bibitem{Burns2020}
K.~J. Burns, G.~M. Vasil, J.~S. Oishi, D.~Lecoanet, B.~P. Brown,
  \href{https://link.aps.org/doi/10.1103/PhysRevResearch.2.023068}{{Dedalus: A
  flexible framework for numerical simulations with spectral methods}}, Phys.
  Rev. Res. 2 (2020) 023068.
\newblock \href {https://doi.org/10.1103/PhysRevResearch.2.023068}
  {\path{doi:10.1103/PhysRevResearch.2.023068}}.
\newline\urlprefix\url{https://link.aps.org/doi/10.1103/PhysRevResearch.2.023068}

\bibitem{Dumbgen2016}
L.~D{\"u}mbgen, Einf{\"u}hrung in die Statistik, Springer, 2016.

\bibitem{tvalues}
\href{https://www.statistik.tu-dortmund.de/fileadmin/user_upload/Lehrstuehle/Oekonometrie/Lehre/WiSoOekoSS17/tabelletV.pdf}{Tabelle
  der t-verteilung}, accessed 26 March 2024 (2017).
\newline\urlprefix\url{https://www.statistik.tu-dortmund.de/fileadmin/user_upload/Lehrstuehle/Oekonometrie/Lehre/WiSoOekoSS17/tabelletV.pdf}

\bibitem{Troeltzsch2010}
F.~Tr{\"o}ltzsch, {Optimal control of partial differential equations: Theory,
  methods, and applications}, Vol. 112, American Mathematical Soc., 2010.

\bibitem{Armijo1966}
L.~Armijo, Minimization of functions having lipschitz continuous first partial
  derivatives., Pacific Journal of Mathematics 16~(1) (1966) 1 – 3.

\bibitem{Kaipio2007}
J.~Kaipio, E.~Somersalo,
  \href{https://www.sciencedirect.com/science/article/pii/S0377042705007296}{Statistical
  inverse problems: Discretization, model reduction and inverse crimes},
  Journal of Computational and Applied Mathematics 198~(2) (2007) 493--504,
  special Issue: Applied Computational Inverse Problems.
\newblock \href {https://doi.org/https://doi.org/10.1016/j.cam.2005.09.027}
  {\path{doi:https://doi.org/10.1016/j.cam.2005.09.027}}.
\newline\urlprefix\url{https://www.sciencedirect.com/science/article/pii/S0377042705007296}

\end{thebibliography}


%
%
%
\end{document}